\documentclass[11pt]{article}
\usepackage{tikz-cd}
\usetikzlibrary{matrix,arrows,decorations.pathmorphing}
\usepackage{amsthm}
\usepackage{amssymb}
\usepackage{amsmath}
\usepackage{amsmath,amscd}
\usepackage[all,cmtip]{xy}

\begin{document}
\title{Derived functor cohomology groups with yoneda Product}

\maketitle

\author{Hafiz Syed Husain\\
Department of Mathematical Sciences,\\
Federal Urdu University of Arts, Science \& Technology, Karachi.\\
Email: hsyed.hussain@fuuast.edu.pk\\
\& \\
Mariam Sultana\\
Department of Mathematical Sciences,\\
Federal Urdu University of Arts, Science \& Technology, Karachi.\\
Email: marium.sultana@fuuast.edu.pk\\

\begin{abstract}

This work presents an exposition of both the internal structure of derived category of an abelian category $D^{*}(\mathcal{A})$ and its contribution in solving problems, particularly in algebraic geometry. Calculation of some morphisms will be presented between objects in $D^{*}(\mathcal{A})$ as elements in appropriate cohomology groups along with their compositions with the help of Yoneda construction under the assumption that the homological dimension of $D^{*}(\mathcal{A})$ is greater than or equal to $2$. These computational settings will then be considered under sheaf cohomological context with a particular case from projective geometry.
\end{abstract}

\textbf{Keywords:} Derived Category, Triangulated Category, Yoneda Product, Sheaf Cohomology, Smooth Projective Variety.\newline \newline
\textbf{Subject Classification:} 13D09, 14F05, 18E30.\newline \newline

\newtheorem{Hus1}{Definition}
\newtheorem{Hus2}{Proposition}
\newtheorem{Hus3}{Lemma}
\newtheorem{Hus4}{Corollary}
\newtheorem{Hus5}{Results and Discussion}
\newcommand{\Hom}{\operatorname{Hom}}
\newcommand{\id}{\operatorname{id}}
\newcommand{\ch}{\operatorname{ch}}
\newcommand{\td}{\operatorname{td}}
\newcommand{\Hilb}{\operatorname{Hilb}}
\newcommand{\Vect}{\operatorname{Vect}}
\newcommand{\Ext}{\operatorname{Ext}}
\newcommand{\Aut}{\operatorname{Aut}}
\newcommand{\Pic}{\operatorname{Pic}}

\section{Introduction}
The introduction to the notion of both, the triangulated and derived categories can be motivated from classical accounts of \cite{bics, rhrd}, or from a comparatively recent and comprehensive exposition \cite{cwha, sigyimmha}. However, for those coming from non-specialist backgrounds, a very concise introduction can be found in \cite{dhfmtag, ordercoh, rtdcwm, tbdercoh} and the appendix of \cite{ubfmtnt}. This roughly amounts to conceiving derived categories $D(\mathcal{A})$, modelled on an abelian category $\mathcal{A}$, as the objects that solve the following universal mapping problem
\[
\begin{tikzcd}[column sep=small]
\text{Kom}(\mathcal{A}) \arrow{r}{Q}  \arrow{rd}{F}
  & D(\mathcal{A}) \arrow{d}{\exists !G} \\
    & D
\end{tikzcd}
\]
where Kom$(\mathcal{A})$ is the category of complexes of objects from $\mathcal{A}$, $Q$ is the functor that maps quasi-isomorphisms to isomorphisms (i.e. whenever $f:A^{\bullet}\rightarrow B^{\bullet}$ a morphism in Kom$(\mathcal{A})$ with $H^{n}(f):H^{n}(A^{\bullet})\stackrel{\simeq}{\rightarrow} H^{n}(B^{\bullet})$, then $Q(f)$ is an invertible arrow in $\mathcal{M}or(D(\mathcal{A}))$ -the collection of all morphisms in $D(\mathcal{A})$); then any functor $F$ from Kom$(\mathcal{A})$ to $D$ that maps quasi-isomorphisms in Kom$(\mathcal{A})$ to isomorphisms in $D$ must factor through $G$ uniquely. This makes derived category merely an object that exists. It is important that one must have a way of actually carrying out the construction that would satisfy the corresponding universal diagram. This is worked out with the help of the notion of localizing class $S$, for which the above functor $Q$ is the localizing functor, a process that pretty much mimic the process of fractionalizing a ring $R$ over a multiplicative set $S$, specifically when the canonical map $\phi:R\rightarrow R[S^{-1}]$ may not be injective.  This whole process admits an interpretation of $D(\mathcal{A})$ as the localization of homotopy category of the same abelian category denoted by $\mathcal{K}(\mathcal{A})$, by class of all quasi-isomorphisms (see \cite{rtdcwm}, \cite{rhrd} Ch. I and \cite{sigyimmha}, III for detail). Once derived category is identified this way, the internal structure of it renders its object be interpreted as complexes of objects from $\mathcal{A}$ but the class of morphisms -denoted by $\mathcal{M}or(D(\mathcal{A}))$ becomes quite complicated to describe as a result of the collapsing that occurs due to the equivalence relation induced from the localizing. This is then captured through the calculus of fractions as both right $S$-roofs and left $S$-roofs (\cite{bics} XI). Here we may assume the facts that every derived category is triangulated and the abelian category $\mathcal{A}$ sits inside the corresponding $D(\mathcal{A})$ as a full subcategory (cf. \cite{rhrd} I). Also, it is usually $D^{*}(\mathcal{A})$ with $* \in \{+, -, b\}$ that has practical significance in most applications corresponding to the interpretation that the derived category in context (up to isomorphism in $D(\mathcal{A})$) consists of objects which are bounded from left, or right, or both from left and right respectively. However, we will be primarily assuming $D^{+}(\mathcal{A})$ unless stated otherwise after the fact that our category $\mathcal{A}$ has enough injectives (however in case $\mathcal{A}$ corresponds to the category of coherent sheaves neither the existence of enough projectives nor enough injectives can be taken for granted; \cite{dhfmtag} II and Appendix \cite{ubfmtnt}). In order to exemplify the internal structure of $D(\mathcal{A})$ what we want to do is calculate and describe morphisms between its objects which are of particular interest in application in geometry. For this we require Yoneda's construction the detail of which can be found in \cite{dhfmtag} $3.3$, \cite{rhrd} I, \cite{ubfmtnt} A and \cite{sigyimmha} III, which predominantly describe the construction over coherent sheaves. This helps interpret its extension to  the derived functors $\text{Ext}^{i}_{\mathcal{A}}(A^{\bullet},B^{\bullet})$ as $\text{Hom}_{D(\mathcal{A})}(A^{\bullet},B^{\bullet}[i])$ defining a grading on Ext groups as $\bigoplus_{i\geq 0}^{n}\text{Ext}_{\mathcal{A}}^{i}(A^{\bullet},B^{\bullet})$ such that $A^{\bullet}, B^{\bullet}, C^{\bullet}$ are successive elements in the (possibly) infinite sequence $\in \Pi_{i\geq n}D^{\bullet}_{i},\; i\in \mathbb{Z}$ of objects from $D(\mathcal{A})$; with Yoneda products
\[
\text{Ext}_{\mathcal{A}}^{i}(A^{\bullet},B^{\bullet})\times \text{Ext}_{\mathcal{A}}^{j}(B^{\bullet},C^{\bullet}) \longrightarrow \text{Ext}_{\mathcal{A}}^{i+j}(A^{\bullet},C^{\bullet})
\]
such that these products are then found to be coinciding with compositions in the derived category $D(\mathcal{A})$

\begin{align*}
&\text{Hom}_{D(\mathcal{A})}(A^{\bullet}[k],B^{\bullet}[k+i])\times \text{Hom}_{D(\mathcal{A})}(B^{\bullet}[k+i],C^{\bullet}[k+i+j])\\
& \longrightarrow \text{Hom}_{D(\mathcal{A})}(A^{\bullet}[k],C^{\bullet}[k+i+j])
\end{align*}
\section{Results and Discussion}
In what follows, we work out some examples of these constructions and interpretations and give our explicit calculations. Here we assume the embedding of $\mathcal{A}$ in $D(\mathcal{A})$ as a full subcategory and apply both $\text{Ext}_{\mathcal{A}}(-,-)$ and $\text{Hom}_{D(\mathcal{A})}(-,-)$ on objects and morphisms from $\mathcal{A}$ considered inside Ab (i.e. the category of abelian groups) and $D(\mathcal{A})$ respectively. This will have the advantage of avoiding the use of spectral sequences which are an indispensable tool when such examples are computed under the identification of Ext and Hom at $D(\mathcal{A})$ in general.
\begin{Hus3}
Let $\mathcal{A}$ be any abelian category with a filtration $\mathcal{F}_{1}\subset \mathcal{F}_{2} \subset \mathcal{G}$ of $\mathcal{G}$ from $\mathcal{A}$ yielding an exact sequence
\begin{equation}\label{1}
0\longrightarrow \mathcal{F}_{1}\longrightarrow \mathcal{F}_{2} \stackrel{f}{\longrightarrow} \mathcal{G} / \mathcal{F}_{1} \longrightarrow \mathcal{G}/\mathcal{F}_{2}\longrightarrow 0
\end{equation}
where $f$ is the usual composition of the embedding $\mathcal{F}_{2}\rightarrow \mathcal{G}$ followed by canonical surjection $\mathcal{G}\rightarrow \mathcal{G}/\mathcal{F}_{1}$. Let $\alpha \in \text{Ext}_{\mathcal{A}}^{2}(\mathcal{G}/\mathcal{F}_{2},\mathcal{F}_{1})\simeq \text{Hom}_{D(\mathcal{A})}(\mathcal{G}/\mathcal{F}_{2},\mathcal{F}_{1}[2])$ correspond to this exact sequence then $\exists \; \;  \alpha_{1} \in \text{Ext}_{\mathcal{A}}^{1}(\mathcal{F}_{2}/\mathcal{F}_{1},\mathcal{F}_{1}) \simeq \text{Hom}_{D(\mathcal{A})}(\mathcal{F}_{2}/\mathcal{F}_{1},\mathcal{F}_{1}[1])$ and $\alpha_{2} \in \text{Ext}_{\mathcal{A}}^{1}(\mathcal{G}/\mathcal{F}_{2},\mathcal{F}_{2}/\mathcal{F}_{1}) \simeq \text{Hom}_{D(\mathcal{A})}(\mathcal{G}/\mathcal{F}_{2},(\mathcal{F}_{2}/\mathcal{F}_{1})[1])$ corresponding to short exact sequences $0\longrightarrow\mathcal{F}_{1}\longrightarrow\mathcal{F}_{2}\longrightarrow \mathcal{F}_{2}/\mathcal{F}_{1}\longrightarrow 0$ and $0\longrightarrow \mathcal{F}_{2}/\mathcal{F}_{1}\longrightarrow \mathcal{G}/\mathcal{F}_{1}\longrightarrow \mathcal{G}/\mathcal{F}_{2}\longrightarrow 0$ respectively; such that $\alpha = \alpha_{1}\cdot \alpha_{2}$

\end{Hus3}
\begin{proof}
Following Yoneda's construction, we first describe $\alpha_{1}$ and $\alpha_{2}$ in $\text{Hom}_{D(\mathcal{A})}(\mathcal{F}_{2}/\mathcal{F}_{1},\mathcal{F}_{1}[1])$ and $\text{Hom}_{D(\mathcal{A})}(\mathcal{G}/\mathcal{F}_{2},\mathcal{F}_{2}/\mathcal{F}_{1}[1])$ respectively. Define an acyclic complex
\[
\mathcal{L}^{\bullet}: 0\rightarrow \mathcal{L}^{-1}=\mathcal{F}_{1}\stackrel{\text{d}^{-1}}{\rightarrow} \mathcal{L}^{0}=\mathcal{F}_{2}\stackrel{\text{d}^{0}}{\rightarrow} \mathcal{L}^{1}=\mathcal{F}_{2}/\mathcal{F}_{1} \rightarrow 0
\]
with $\text{d}^{-1}$ corresponds to the usual embedding and $\text{d}^{0}$ corresponds to the canonical projection, then we obtain a quasi-isomorphism: $\tilde{\mathcal{L}}^{\bullet}\stackrel{s_{\mathcal{L}^{\bullet}}}{\longrightarrow} (\mathcal{F}_{2}/\mathcal{F}_{1})[0]$, such that
\[
\tilde{\mathcal{L}}^{\bullet}: 0\longrightarrow \mathcal{L}^{-1}=\mathcal{F}_{1} \stackrel{\text{d}^{-1}}{\longrightarrow} \mathcal{L}^{0}=\mathcal{F}_{2} \longrightarrow 0
\]
and $s_{\mathcal{L}^{\bullet}}^{i}=0, \; \forall i\neq 0, s_{\mathcal{L}^{\bullet}}^{i}=\text{d}^{0}_{\mathcal{L}^{\bullet}}, i=0$
similarly we get $\tilde{\mathcal{L}}^{\bullet}\stackrel{g_{\mathcal{L}^{\bullet}}}{\longrightarrow} \mathcal{F}_{1}[1]$ such that $g_{\mathcal{L}^{\bullet}}^{i}=0,\; \forall i\neq -1, g_{\mathcal{L}^{\bullet}}^{-1}=\text{id}_{\mathcal{F}_{1}}$.

 This gives both $\alpha_{1}$ and $T(\alpha_{1})=\alpha_{1}[1]$ as left S-roofs (or as a left fraction $(s_{\mathcal{L^{\bullet}}},g_{\mathcal{L^{\bullet}}})=:g_{\mathcal{L^{\bullet}}}/s_{\mathcal{L^{\bullet}}}$ ), with $T$ being the usual autoequivalence of shift functor that ${D(\mathcal{A})}$ inherits from its triangulated structure (thus $\text{Hom}_{D(\mathcal{A})}(A^{\bullet},B^{\bullet})\simeq \text{Hom}_{D(\mathcal{A})}(T(A^{\bullet}),T(B^{\bullet}))$ and $A^{\bullet}[1]\simeq T(A^{\bullet})$ as an object in $D(\mathcal{A})$)
\[
\alpha_{1}:\; (\mathcal{F}_{2}/\mathcal{F}_{1})[0]\stackrel{s_{\mathcal{L^{\bullet}}}}{\longleftarrow} \tilde{\mathcal{L}}^{\bullet}\stackrel{g_{\mathcal{L^{\bullet}}}}{\longrightarrow} \mathcal{F}_{1}[1]\Longrightarrow T(\alpha_{1}):\; (\mathcal{F}_{2}/\mathcal{F}_{1})[1]\stackrel{s_{\mathcal{L^{\bullet}}}[1]}{\longleftarrow} \tilde{\mathcal{L}}^{\bullet}[1]\stackrel{g_{\mathcal{L^{\bullet}}}[1]}{\longrightarrow} \mathcal{F}_{1}[2]
\]
Now consider $\alpha_{2} \in \text{Hom}_{D(\mathcal{A})}(\mathcal{G}/\mathcal{F}_{2},\mathcal{F}_{2}/\mathcal{F}_{1}[1])$. We similarly define
\[
\mathcal{K}^{\bullet}: 0\rightarrow \mathcal{K}^{-1}=\mathcal{F}_{2}/\mathcal{F}_{1} \stackrel{\text{d}^{-1}}{\rightarrow} \mathcal{K}^{0}=\mathcal{G}/\mathcal{F}_{1} \stackrel{\text{d}^{0}}{\rightarrow} \mathcal{K}^{1}=\mathcal{G}/\mathcal{F}_{2} \rightarrow 0
\]
where $\text{d}^{-1}$ and $\text{d}^{0}$ are canonical projections which are guaranteed from the third isomorphism theorem of modules over a ring once Fred-Mitchell Embedding theorem is assumed (see \cite{cwha} $1.6.1$); giving us the quasi-isomorphism: $\tilde{\mathcal{K}}^{\bullet}\stackrel{s_{\mathcal{K}^{\bullet}}}{\longrightarrow} (\mathcal{G}/\mathcal{F}_{2})[0] $, with $s_{\mathcal{K}^{\bullet}}^{i}=0, \forall i\neq 0, s_{\mathcal{K}^{\bullet}}^{0}=\text{d}^{0}_{\mathcal{K}^{\bullet}}$ and $\tilde{\mathcal{K}}^{\bullet}\stackrel{g_{\mathcal{K}^{\bullet}}}{\longrightarrow}(\mathcal{F}_{2}/\mathcal{F}_{1})[1]$ given by $g_{\mathcal{K}^{\bullet}}^{i}=0, \;\forall i\neq -1, g_{\mathcal{K}^{\bullet}}^{-1}=\text{id}_{(\mathcal{F}_{2}/\mathcal{F}_{1})[1]} $.
This gives us the representation of $\alpha_{2}$ as a left S-roof (or as a left fraction $(s_{\mathcal{K^{\bullet}}},g_{\mathcal{K^{\bullet}}})=:g_{\mathcal{K^{\bullet}}}/s_{\mathcal{K^{\bullet}}}$ ) as follows
\[
(\mathcal{G}/\mathcal{F}_{2})[0]\stackrel{s_{\mathcal{K^{\bullet}}}}{\longleftarrow} \tilde{\mathcal{K}}^{\bullet}\stackrel{g_{\mathcal{K^{\bullet}}}}{\longrightarrow} (\mathcal{F}_{2}/\mathcal{F}_{1})[1]
\]
We can now define the composition $\alpha=\alpha_{1}\cdot \alpha_{2}=y(\mathcal{L^{\bullet}})\circ y(\mathcal{K^{\bullet}}) \Rightarrow (s_{\mathcal{L^{\bullet}}}[1],g_{\mathcal{L^{\bullet}}}[1])\circ (s_{\mathcal{K^{\bullet}}},g_{\mathcal{K^{\bullet}}})$ such that we get

\[
\xymatrix{
 &  &\tilde{\mathcal{M}}^{\bullet} \ar[ld]_{s^{\prime}} \ar[rd]^{g^{\prime}} &  & \\
  &\tilde{\mathcal{K}}^{\bullet}\ar[ld]_{s_{\mathcal{K^{\bullet}}}} \ar[rd]^{g_{\mathcal{K^{\bullet}}}}  &  &\tilde{\mathcal{L}}^{\bullet}[1]\ar[ld]_{s_{\mathcal{L^{\bullet}}}[1]} \ar[rd]^{g_{\mathcal{L^{\bullet}}}[1]}  &\\
(\mathcal{G}/\mathcal{F}_{2})[0] &  &(\mathcal{F}_{2}/\mathcal{F}_{1})[1]  &  &\mathcal{F}_{1}[2]}
\]
with $\tilde{\mathcal{M}}^{\bullet}=y(\mathcal{M})$ such that $\mathcal{M}$ is the same length $2$ exact sequence given in $(1)$ above, giving $\alpha=\alpha_{1}\cdot \alpha_{2}=(s_{\mathcal{K^{\bullet}}}s^{\prime}, g_{\mathcal{L^{\bullet}}}[1]g^{\prime}) \in \text{Hom}_{D(\mathcal{A})}(\mathcal{G}/\mathcal{F}_{2},\mathcal{F}_{1}[2])$

\end{proof}

\begin{Hus3}
Let $\alpha \in \text{Ext}_{\mathcal{A}}^{2}(\mathcal{G}/\mathcal{F}_{2},\mathcal{F}_{1})$  be as proposed in Lemma $1$ above, then $\alpha = 0$.
\end{Hus3}

\begin{proof}
We show $\alpha = 0$ by showing that $\alpha \simeq 0$ in $\mathcal{M}or(D(\mathcal{A}))$. Now from Proposition $1$ above we know that $\alpha$ in $\text{Hom}_{D(\mathcal{A})}(\mathcal{G}/\mathcal{F}_{2}, \mathcal{F}_{1}[2])\subset \mathcal{M}or(D(\mathcal{A}))$ is representable as a left S-roof
\[
\xymatrix{
   &\tilde{\mathcal{M}}^{\bullet} \ar[ld]_{s_{\mathcal{K}^{\bullet}}s^{\prime}} \ar[rd]^{g_{\mathcal{L^{\bullet}}}g^{\prime}} &   \\
  (\mathcal{G}/\mathcal{F}_{2})[0] &  &\mathcal{F}_{1}[2]}
\]
We define a morphism $0\in \text{Hom}_{D(\mathcal{A})}(\mathcal{G}/\mathcal{F}_{2}, \mathcal{F}_{1}[2])$ as a left S-roof
\[
\xymatrix{
   &\mathcal{J}^{\bullet} \ar[ld]_{\text{id}} \ar[rd]^0 &   \\
  (\mathcal{G}/\mathcal{F}_{2})[0] &  &\mathcal{F}_{1}[2]}
\]
with $(\mathcal{J}^{\bullet})=(\mathcal{G}/\mathcal{F}_{2})[0]$, then we can have $\mathcal{I}^{\bullet} : 0\rightarrow 0\rightarrow \mathcal{F}_{2} \rightarrow \mathcal{G}/\mathcal{F}_{1} \rightarrow 0$ such that we get
\[
\xymatrix{
   &\mathcal{I}^{\bullet} \ar[ld]_{s_{\mathcal{I}^{\bullet}}} \ar[rd]^{g_{\mathcal{I}^{\bullet}}} &   \\
  \tilde{\mathcal{M}}^{\bullet} &  &\mathcal{J}^{\bullet} }
\]
with $s_{\mathcal{I}^{\bullet}}^{i} = \text{id}$ for $0\leq i\leq 1$ and is trivial in all other $i$, $g_{\mathcal{I}^{\bullet}}^{i}=\text{id}$ for $i=0$ and is trivial for all other $i$. All this finally gives rise to a following diagram that commutes up to homotopy (which is fairly straightforward to verify)
\[
\xymatrix{
 &(\mathcal{G}/\mathcal{F}_{2})[0]  &\\
 \tilde{\mathcal{M}}^{\bullet} \ar[ru]^{s_{\mathcal{K}^{\bullet}}s^{\prime}} \ar[rd]_{g_{\mathcal{L^{\bullet}}}g^{\prime}} &\mathcal{I}^{\bullet} \ar[l]_{s_{\mathcal{I}^{\bullet}}} \ar[r]^{g_{\mathcal{I}^{\bullet}}} &\mathcal{J}^{\bullet} \ar[ld]^0 \ar[lu]_{\text{id}}\\
   &\mathcal{F}_{1}[2]  & }
\]
thereby establishing that $\alpha \simeq 0$ in $\text{Hom}_{D(\mathcal{A})}(\mathcal{G}/\mathcal{F}_{2},\mathcal{F}_{1}[2])$, which implies $\alpha = 0$ in $\text{Ext}^{2}(\mathcal{G}/\mathcal{F}_{2},\mathcal{F}_{1})$.
\end{proof}
This gives $\alpha$ as a  $2$-dimensional analogue of $0$ in $\text{Ext}^{1} (-,-)$ which can be interpreted analogously as a short exact sequence that splits. The classical homological approach of Cartan-Eilenberg would have accessed the same situation with the help of resolutions (injective or projective) giving Ext-groups as derived functor cohomology objects $\text{Ext}_{\mathcal{A}}^{i}(C,A)\simeq H^{i}(R(\text{Hom}_{\mathcal{A}}(C,A)))\simeq H^{i}(\text{Hom}_{\mathcal{A}}(C,I^{\bullet}_{A}))$ or  $\text{Ext}_{\mathcal{A}}^{i}(C,A)\simeq H^{i}(R(\text{Hom}_{\mathcal{A}}(C,A)))\simeq H^{i}(\text{Hom}_{\mathcal{A}}(P^{\bullet}_{C},A))$; where $A\stackrel{quis}{\rightarrow} I^{\bullet}$ and $C\stackrel{quis}{\rightarrow} P^{\bullet}$ are injective and projective resolutions of $A$ and $C$ respectively which are quasi-isomorphisms, hence isomorphisms in $D(\mathcal{A})$ ($A, C\in \text{Ob}\mathcal{A}$). Even at this level without historical reference to derived category, one can see the transition from $\text{Kom}(\mathcal{A})$ to $D(\mathcal{A})$ via $\mathcal{K}(\mathcal{A})$ as follows:
consider
\[ \text{Ext}^{i}_{\mathcal{A}}(C,A)= \frac{\text{ker}(\text{Hom}_{\text{Ab}}(C, I^{i})\stackrel{\text{d}_{*I^{\bullet}}^{i}}{\rightarrow}\text{Hom}_{\text{Ab}}(C,I^{i+1}))}{\text{Im}
(\text{Hom}_{\text{Ab}}(C, I^{i-1})\stackrel{\text{d}_{*I^{\bullet}}^{i-1}}{\rightarrow}\text{Hom}_{\text{Ab}}(C,I^{i}))}
\]
then $\alpha \in \text{Ext}_{\mathcal{A}}^{i}(C,A) \Rightarrow \alpha \in \text{Hom}_{\text{Ab}}(C,I^{i})$ determining the morphism of complexes $\alpha^{\bullet}\in \mathcal{M}or(\text{Kom}(\mathcal{A}))$ as $\alpha^{\bullet}: C\rightarrow I^{\bullet}, \alpha^{j}=0, \; \forall j\neq i, \alpha^{i}=\alpha:C\rightarrow I^{i}$. However, if $\alpha=0 \in \text{Ext}_{\mathcal{A}}^{i}(C,A)$, then it determines a morphism of complexes $\alpha^{\bullet}\in \mathcal{M}or(\mathcal{K}(\mathcal{A}))$ such that $\alpha^{\bullet}=0\in \text{Hom}_{\mathcal{K}(\mathcal{A})}(C,I^{\bullet}[i])$, since it is homotopically trivial determined by homotopies $h=\Pi_{j\geq 0} h^{j}, h^{j}: C\rightarrow I^{j-1}$ such that $h^{j}=0, \; \forall j\neq i, h^{i}= \beta \in \text{Hom}_{\text{Ab}}(C,I^{i-1})$ with $\text{d}_{I^{\bullet}}^{i-1}(\beta)=0$. We thus have $\text{Ext}^{i}_{\mathcal{A}}(C,A)\simeq \text{Hom}_{D(\mathcal{A})}(C,A[i])$ via $\text{Ext}^{i}_{\mathcal{A}}(C,A)\simeq \text{Hom}_{\mathcal{K}(\mathcal{A})}(C,I^{\bullet}[i])$, $\text{Hom}_{\mathcal{K}(\mathcal{A})}(C,I^{\bullet}[i])\simeq \text{Hom}_{D(\mathcal{A})}(C,I^{\bullet}[i])$ and $\text{Hom}_{D(\mathcal{A})}(C,I^{\bullet}[i])\simeq \text{Hom}_{D(\mathcal{A})}(C,A[i]))$. All this then further generalizes to genuine complexes via the notion of inner-hom $\text{Hom}^{\bullet}(A^{\bullet}, B^{\bullet})$ which is a complex in Kom(Ab) the $n^{th}$-objecct of which is $(\text{Hom}^{\bullet}(A^{\bullet}, B^{\bullet}))^{n}=\Pi _{i\in \mathbb{Z}}\text{Hom}_{Ab}(A^{i},B^{i+n})$ with differentials $d^{n}(-)=d_{B^{\bullet}}(-)-(-1)^{n}(-)d_{A^{\bullet}}$, which determines the homotopy equivalence to $0$ at each i$^{th}$ degree, the special case of which was already determined above as $\text{d}_{I^{\bullet}}^{i-1}(\beta)=0$. For instance $f\in (\text{Hom}^{\bullet}(A^{\bullet}, B^{\bullet}))^{n} \Rightarrow f=(\ldots, f_{-1}, f_{0}, f_{1},\ldots )$, then at $n^{th}$ place, $d^{n}(f)=d_{B^{\bullet}}^{i+n}(f_{i})-(-1)^{n}(f_{i})d_{A^{\bullet}}^{i-1}$, here $n$ is just the shift functor in disguise that would later descend to triangulated $\mathcal{K}(\mathcal{A})$ and then to $D(\mathcal{A})$ and the negative sign in $(-1)^{n}$ takes care of the negative of the differential in $A^{\bullet}$ which is the result of shifting the complex by $i^{th}$-degree to the left (depending upon whether $i$ is even or odd) to match the homotopy calculation. Then we define $\text{Ext}^{i}(A^{\bullet},B^{\bullet})\simeq H^{i}(R(\text{Hom}^{\bullet}(A^{\bullet}, B^{\bullet})))$ as a convergent spectral sequence (cf. \cite{dhfmtag} II).

We want to discuss a concrete application to match the abstract setting of Lemma $1$ and $2$ from algebraic geometry (which will be our Proposition $2$ below). From this point onwards we define a smooth complex projective algebraic variety $X$ as a separated scheme of finite type over $\mathbb{C}$ which is projective over $\mathbb{C}$  and is locally regular; i.e. all its local rings are regular (we have relaxed the classical Hartshorne condition of integrality \cite{rhag} II.$4$); thus all our varieties $X$ are smooth, complex and projective. Let $\mathcal{A}$ denote the category of coherent sheaves on $X$. Then $\mathcal{A}$ can already be seen as a weak invariant of $X$ in the sense if dim$(X)=n$ then homological dimension of $\mathcal{A}$ -denoted by $dh(\mathcal{A})$ equals $n$ as well and the fact that if there is an equivalence between categories of coherent sheaves of $X_{1}$ and $X_{2}$ then $X_{1}\simeq X_{2}$; a corollary and an application Orlov's famous theorem about the existence and uniqueness of Fourier-Mukai transforms (cf. \cite{ordercoh}). Also, we will not make any distinction between classical notion of a variety and its scheme theoretic counterpart, something which makes perfect sense from \cite{rhag} II, $2.4$. Following proposition has its motivation in the classical sources as \cite{rhrd} and \cite{rhag}. We give our explicit proof as follows:
\begin{Hus2}
Let $X$ be any smooth complex projective variety with $D^{b}(X)$ as its bounded derived category of coherent sheaves and $\omega_{X}$ be its canonical bundle. Let $R\mathcal{H}om(\mathcal{F}^{\bullet}, \mathcal{O}_{X})=:_{\text{df}}(\mathcal{F}^{\bullet})^{\vee}$ denote the derived dual of $(\mathcal{F}^{\bullet})$ and $\Pi_{i\in I}(\mathcal{F}_{i}^{\bullet})$ be any sequence of objects from $D^{b}(X)$ with Yoneda products
\[
y: \text{Ext}^{k}_{\mathcal{A}}(\mathcal{F}_{i-1}^{\bullet}, \mathcal{F}_{i}^{\bullet})\times \text{Ext}^{l}_{\mathcal{A}}(\mathcal{F}_{i}^{\bullet}, \mathcal{F}_{i+1}^{\bullet})\longrightarrow \newline \text{Ext}^{k+l}_{\mathcal{A}}(\mathcal{F}_{i-1}^{\bullet}, \mathcal{F}_{i+1}^{\bullet})
\]
then $y$ descends to corresponding products on sheaf cohomology
\begin{align*}
&\tilde{y}:\text{H}^{n-k}(X,(\mathcal{F}_{i-1}^{\bullet})^{\vee} \otimes \mathcal{F}_{i}^{\bullet}\otimes \omega_{X})\times {H}^{n-l}(X,(\mathcal{F}_{i}^{\bullet})^{\vee} \otimes \mathcal{F}_{i+1}^{\bullet}\otimes \omega_{X})\\
&\longrightarrow {H}^{n-k-l}(X,(\mathcal{F}_{i-1}^{\bullet})^{\vee} \otimes \mathcal{F}_{i+1}^{\bullet}\otimes \omega_{X})
\end{align*}
\end{Hus2}

\begin{proof}
First of all we know that $D^{b}(X)$ is $\mathbb{C}$-linear, thus it is equipped with Serre's functor (cf. \cite{bvbgenrep, ordercoh}) the special case of which is Serre's duality $\text{Ext}^{i}(\mathcal{F}, \omega_{X})\simeq \text{H}^{n-i}(X, \mathcal{F})^{*}$ (\cite{rhag} III,$7$). Also since $dh(D^{b}(X))=\text{dim}(X)$ thus all gradings $\bigoplus_{0\leq i\leq n}\text{Ext}_{\mathcal{A}}^{i}$ are finite, such that $\text{Ext}^{p+q}(A^{\bullet},B^{\bullet})$ is obtained as the convergent spectral sequence $E_{2}^{p,q} = R^{p}\text{Hom}^{\bullet}(A^{\bullet},H^{q}(B)^{\bullet})\simeq \text{Ext}^{p}(A^{\bullet},H^{q}(B)^{\bullet})$ (provided $\mathcal{A}$ has enough injectives). On the other hand since $\mathcal{A}$ almost never has enough projectives thus we have another spectral sequence $H^{p}_{II}H^{q}_{I}(L^{\bullet,\bullet})=\text{Ext}^{p}(H^{-q}(A^{\bullet}),(B)^{\bullet})$ converging to $\text{Ext}^{p+q}(A^{\bullet},B^{\bullet})$; where $L^{\bullet,\bullet}$ is the Cartan-Eilenberg resolution of $A^{\bullet}$ (see \cite{dhfmtag} $2.2$ or we refer classical exposition \cite{cartanha} for detail). Hence the transition from Ext-groups of coherent sheaves sitting as $0$-degree in $D^{b}(X)$ to Ext-groups of complexes of them makes sense. These identifications through spectral sequences (also known as special cases of Grothendieck spectral sequences) establish $\text{Ext}^{i}_{\mathcal{A}}(A^{\bullet},B^{\bullet})\simeq \text{Hom}_{D(\mathcal{A})}(A^{\bullet},B^{\bullet}[i])$. One similarly uses Leray spectral sequence to arrive at $\text{H}^{i}(X,\mathcal{F}^{\bullet})$ from $\text{H}^{i}(X,\mathcal{F})$ to obtain Serre duality at the level of $D^{b}(X)$; i.e. $\text{Ext}^{i}(\mathcal{F}^{\bullet}, \omega_{X})\simeq \text{H}^{n-i}(X, \mathcal{F}^{\bullet})^{*}$; cf. \cite{rhrd} II and V. Serre's duality from \cite{rhag} III,$7$ only makes sense at the level of $\mathcal{A}$; i.e coherent sheaves. Then from \cite{rhag} III,$6$, and above, we get the identifications
\begin{align*}
\text{Ext}^{i}(\mathcal{F}^{\bullet}, \mathcal{G}^{\bullet}) & \simeq \text{Ext}^{n-i}(\mathcal{G}^{\bullet}, \mathcal{F}^{\bullet}\otimes \omega_{X})\\
& \simeq \text{Ext}^{n-i}(\mathcal{O}_{X}, (\mathcal{G}^{\bullet})^{\vee}\otimes \mathcal{F}^{\bullet}\otimes \omega_{X})\\
& \simeq R^{n-i}\Gamma(X,(\mathcal{G}^{\bullet})^{\vee}\otimes \mathcal{F}^{\bullet}\otimes \omega_{X})\\
&\simeq \text{H}^{n-i}(X,(\mathcal{G}^{\bullet})^{\vee}\otimes \mathcal{F}^{\bullet}\otimes \omega_{X})
\end{align*}
(Note: no need to derive tensor product in above expressions because $\omega_{X}$ is a line bundle.)
\end{proof}

We now relate both Lemma $1$ and $2$ in context of Proposition $1$, and discuss a specific concrete case from projective geometry. It is here that the restriction on homological dimension is particularly informative. Restricting homological dimension to be at least greater than or equal to $2$ involves the cases of complex algebraic surfaces and their higher dimensional analogues only -for instance Calabi-Yau n-folds, since the case of algebraic curves would render all second extensions trivial and thus any Yodeda product and its descent to sheaf cohomolgy as wworked out in Proposition $1$ above will always yield trivial results; a consequence of Grothendieck vanishing theorem (cf. \cite{rhag} III, 2.7). Assuming \cite{rhag} $2.6$, we will not be making any distinction between classical projective $n$-space $\mathbb{P}^{n}_{K}$ over an algebraically closed field $K$ and its scheme theoretic interpretation,  which comes equipped with scheme morphism $\text{Proj}(K[x_{0},x_{2},\ldots ,x_{n}])\longrightarrow \text{Spec}(K)$. The advantage of the latter is that it helps make use of sheaf theoretic tools for (co)-homological computations. We fix $n=3$ and assume $K=\mathbb{C}$ (although much of what this paper presents is valid for any $K$ algebraically closed). Let $\Omega_{\mathbb{P}^{3}_{\mathbb{C}}/\mathbb{C}}$ denote the sheaf of differentials on $\mathbb{P}^{3}_{\mathbb{C}}$, $\mathcal{O}_{\mathbb{P}^{3}_{\mathbb{C}}}(m)$ the structure sheaf of $\mathbb{P}^{3}_{\mathbb{C}}$ twisted $m$-times,$\forall m \in \mathbb{Z}$, by the twisting sheaf of Serre and $\mathcal{T}_{\mathbb{P}^{3}_{\mathbb{C}}}$ the corresponding tangent sheaf. Let $Y$ be a quadric surface in $\mathbb{P}^{3}_{\mathbb{C}}$, then we know that $Y$ can be realized as an image of Segre embedding $i:\mathbb{P}^{1}_{\mathbb{C}}\times \mathbb{P}^{1}_{\mathbb{C}}\longrightarrow \mathbb{P}^{3}_{\mathbb{C}}$ (cf.\cite{jhfcag} $2.7$ p. $25$), with $i_{*}$ denoting the push-forward functor from $\text{Coh}(\mathbb{P}^{1}_{\mathbb{C}}\times \mathbb{P}^{1}_{\mathbb{C}})$ to $\text{Coh}(\mathbb{P}^{3}_{\mathbb{C}})$. On the other hand $Y$ considered as a divisor in divisor class group of $\mathbb{P}^{3}_{\mathbb{C}}$ has its associated line bundle representation, say $\mathcal{L}(D_{Y})\in \text{Pic}(\mathbb{P}^{3}_{\mathbb{C}})$ in Picard group of $\mathbb{P}^{3}_{\mathbb{C}}$, where $D_{Y} \in \text{Cl}(\mathbb{P}^{3}_{\mathbb{C}})$ denotes the corresponding divisor to $Y$ in the divisor class group of $\mathbb{P}^{3}_{\mathbb{C}}$ which is parametrized by $\mathbb{Z}$ such that  $\mathcal{L}(D_{Y})\simeq \mathcal{O}_{\mathbb{P}^{3}_{\mathbb{C}}}(-2)$ (\cite{rhag} $2.6$). For notational brevity, we drop the subscript for $\mathbb{C}$. Then we have the following:

 \begin{Hus2} Let $\tilde{\alpha_{1}} \in \text{H}^{2}(\mathbb{P}^{3},i_{*}(\mathcal{O}_{\mathbb{P}^{1}\times \mathbb{P}^{1}})(-6))$ be the cocycle corresponding to $\alpha_{1} \in \text{Ext}^{1}(i_{*}(\mathcal{O}_{\mathbb{P}^{1}\times \mathbb{P}^{1}}),\mathcal{O}_{\mathbb{P}^{3}}(-2))$ determined by the short exact sequence
 \[
 \alpha_{1}: 0\longrightarrow \mathcal{O}_{\mathbb{P}^{3}}(-2)\longrightarrow \mathcal{O}_{\mathbb{P}^{3}}\longrightarrow i_{*}(\mathcal{O}_{\mathbb{P}^{1}\times \mathbb{P}^{1}})\longrightarrow 0
 \]
 and let $\tilde{\alpha_{2}} \in \text{H}^{2}(\mathbb{P}^{3},(i_{*}(i^{*}(\Omega_{\mathbb{P}^{3}}(-4))))^{\vee})$ be the cocycle corresponding to $\alpha_{2} \in \text{Ext}^{1}(\mathcal{T}_{\mathbb{P}^{3}},i_{*}(\mathcal{O}_{\mathbb{P}^{1}\times \mathbb{P}^{1}}))$ determimned by the short exact sequence
 \[
 \alpha_{2}: 0\longrightarrow i_{*}(\mathcal{O}_{\mathbb{P}^{1}\times \mathbb{P}^{1}})\longrightarrow (\bigoplus_{i=1}^{i=4}\mathcal{O}_{\mathbb{P}^{3}}(1))/\mathcal{O}_{\mathbb{P}^{3}}(-2) \longrightarrow (\bigoplus_{i=1}^{i=4}\mathcal{O}_{\mathbb{P}^{3}}(1))/\mathcal{O}_{\mathbb{P}^{3}}\longrightarrow 0
 \]
 then the Yoneda product as descended to sheaf cohomology in Proposition $1$ above yields
\[
\tilde{y}:\text{H}^{2}(\mathbb{P}^{3},i_{*}(\mathcal{O}_{\mathbb{P}^{1}\times \mathbb{P}^{1}})(-6))\times \text{H}^{2}(\mathbb{P}^{3},(i_{*}(i^{*}(\Omega_{\mathbb{P}^{3}}(-4)))^{\vee})
\longrightarrow \text{H}^{1}(\mathbb{P}^{3},\Omega_{\mathbb{P}^{3}}(-5))
\]
such that the corresponding cohomology product $\tilde{\alpha_{1}}\tilde{y}\tilde{\alpha_{2}}=P(x_{0},\ldots ,x_{3})\tilde{\alpha_{2}}=0 \in \text{H}^{2}(\mathbb{P}^{3},\Omega_{\mathbb{P}^{3}}(-5))$, amounts to multiplication by a degree $10$ monomial $P(x_{0},\ldots ,x_{3})$ to produce a trivial cocycle in ${H}^{1}(\mathbb{P}^{3},\Omega_{\mathbb{P}^{3}}(-5))$.
\end{Hus2}

\begin{proof}
Given that $D_{Y}\in \text{Cl}(\mathbb{P}^{3})$ corresponds to $Y$ as a degree $2$ divisor, it canonically yields the short exact sequence $\alpha_{1}$, whereas $\alpha_{2}$ is obtained by dualizing Euler-sequence corresponding to the sheaf of differentials on $\mathbb{P}^{3}$ and the filtration $\mathcal{O}_{\mathbb{P}^{3}}(-2)\subsetneq \mathcal{O}_{\mathbb{P}^{3}} \subsetneq (\bigoplus_{i=1}^{i=4}\mathcal{O}_{\mathbb{P}^{3}}(1))$ (cf. \cite{rhag} II, 8). Also from this fitlration and the fact that coherent sheaves are locally nothing but finitely generated modules one can glue the local data at the level of any chosen trivialization of, say $\mathbb{P}^{3}$, one obtains the length $2$ sequence as in Lamma 1. Then, since none of these sequences split, both $\text{Ext}^{1}(i_{*}(\mathcal{O}_{\mathbb{P}^{1}\times \mathbb{P}^{1}}),\mathcal{O}_{\mathbb{P}^{3}}(-2))\neq 0$ and $\text{Ext}^{1}(\mathcal{T}_{\mathbb{P}^{3}},i_{*}(\mathcal{O}_{\mathbb{P}^{1}\times \mathbb{P}^{1}}))\neq 0 $; we obtain

\begin{align*}
\text{Ext}^{1}(i_{*}(\mathcal{O}_{\mathbb{P}^{1}\times \mathbb{P}^{1}}),\mathcal{O}_{\mathbb{P}^{3}}(-2)) & \simeq \text{Ext}^{2}(\mathcal{O}_{\mathbb{P}^{3}},i_{*}(\mathcal{O}_{\mathbb{P}^{1}\times \mathbb{P}^{1}})(-6))\\
& \simeq \text{H}^{2}(\mathbb{P}^{3},i_{*}(\mathcal{O}_{\mathbb{P}^{1}\times \mathbb{P}^{1}})(-6))\\
& \simeq \text{H}^{2}(\mathbb{P}^{1}\times \mathbb{P}^{1},\mathcal{O}_{\mathbb{P}^{1}\times \mathbb{P}^{1}}(-6))\\
&\simeq \text{H}^{0}(\mathbb{P}^{1}\times \mathbb{P}^{1},\mathcal{O}_{\mathbb{P}^{1}\times \mathbb{P}^{1}}(3))\\
& \simeq S_{10} \qquad \qquad \qquad \qquad \qquad \qquad \qquad \qquad (1).
\end{align*}
where 3rd isomorphism is due to a corollary of Leray Spectral Sequence (\cite{dhfmtag} p74) and the rest are results from standard cohomology of projective space such that $S=\mathbb{C}[x_{0},\ldots x_{3}]$ is the graded ring with $S_{10}$ corresponding to the degree $10$ part (cf. \cite{rhag} 5.1). One similarly obtains

\begin{align*}
\text{Ext}^{1}(\mathcal{T}_{\mathbb{P}^{3}}, i_{*}(\mathcal{O}_{\mathbb{P}^{1}\times \mathbb{P}^{1}}) & \simeq \text{Ext}^{2}(\mathcal{O}_{\mathbb{P}^{3}},\mathcal{T}_{\mathbb{P}^{3}}\otimes i_{*}(\mathcal{O}_{\mathbb{P}^{1}\times \mathbb{P}^{1}})^{\vee}(-4))\\
& \simeq \text{H}^{1}(\mathbb{P}^{3},(\Omega_{\mathbb{P}^{3}}(-4)\otimes i_{*}(\mathcal{O}_{\mathbb{P}^{1}\times \mathbb{P}^{1}}))^{\vee})\\
& \simeq \text{H}^{1}(\mathbb{P}^{3},i_{*}(i^{*}(\Omega_{\mathbb{P}^{3}}(-4)))^{\vee}) \qquad \qquad \qquad \qquad (2).
\end{align*}

All of the above isomorphisms are standard algebraic geometric manipulations (cf. \cite{rhag} II,1, III, 5 and 6 for detail). Since double dual of any algebraic vector bundle is the bundle itself, one obtains $\text{Ext}^{2}(\mathcal{T}_{\mathbb{P}^{3}}, \mathcal{O}_{\mathbb{P}^{3}}(-2))\simeq \text{H}^{1}(\mathbb{P}^{3},\Omega_{\mathbb{P}^{3}}(-5))$ . Thus applying Lemma $1$, $2$ and Proposition $1$ on $\alpha_{1}$ and $\alpha_{2}$ and combining isomorphism (1) and (2) above one obtains the result.

\end{proof}
\subsection{Conclusion}
Lemma $1$ and $2$ present insight about internal structure of derived categories $D^{*}(A)$ in general. This is comparatively the algebraic side of this work. Proposition $1$ presents the descent of Yoneda product from Ext groups to sheaf cohomology, which may correspondingly be considered the more geometric side. Proposition $2$ then connects all of these results in a very concrete projective geometric setting. This investigation thus leads to an insight of how homological methods from more algebraic settings are applicable to geometry.

\section{Statement of Conflict of Interest}
It is hereby declared that there is no conflict of interest between the author and any third party.

\end{document}